%% file: main.tex
\documentclass[Afour,sageh,times,doublespace]{sagej}

\input{preamble}

\begin{document}

\runninghead{Lee et. al.}

\title{ORCHA: A Performance Portability System for Extreme Heterogeneity}

\author {Youngjun Lee\affilnum{1}, Klaus Weide\affilnum{1,2}, Wesley
  Kwiecinski\affilnum{3}, Jared O'Neal\affilnum{1}, Johann
  Rudi\affilnum{4} and Anshu Dubey\affilnum{1,2}}
\affiliation{\affilnum{1} Argonne National Laboratory, USA \\
  \affilnum{2} The University of Chicago, USA \\
  \affilnum{3} University of Illinois Chicago, USA \\
  \affilnum{4} Virginia Tech, USA }

\corrauth{Youngjun Lee, Mathematics and Computer Science Division
  Argonne National Laboratory,
  9700 S. Cass Ave, Lemont IL 60439, USA}
\email{leey@anl.gov}

\begin{abstract}
Heterogeneity is the prevalent trend in the rapidly evolving
high-performance computing (HPC) landscape in both hardware and
application software. The diversity in hardware platforms, currently comprising
various accelerators and a future possibility of specializable chiplets,
poses a significant challenge for scientific software developers
aiming to harness optimal performance across different computing
platforms while maintaining the quality of solutions when their
applications are simultaneously growing more complex. Code synthesis
and code generation can provide mechanisms to mitigate this 
challenge. We have developed
a divide and conquer approach where different aspects of performance are handled by different stand-alone tools that are interfaced with the application through generated code. This portability system, \orcha, enables users to configure and orchestrate their computations among available resources on a platform by specifying a high-level recipe, thereby permitting a {\em many-to-many} paradigm where each recipe results in a different variant of the application.
The core design goal is to let users decide the application's hardware mapping and orchestration by editing only the high-level recipe—without modifying the maintained source code or binding the application to a particular runtime system.
Tools in \orcha distribution are: \cgkit for translating the recipe into an execution graph; \milhoja to execute the graph by orchestrating data and task movement among hardware resources; and \macro that enables users to define their own code-shorthand for higher composability and easier management of code variants. 
Additionally, the design of \orcha permits tools to work in a plug-and-play mode where the application can build and run without \cgkit and \milhoja, and either tool can be swapped out for other tools with similar capabilities by modifying the code generation portion of \orcha. In this paper, we describe the design of \orcha and the role that code-generation plays in isolating applications from tools. We demonstrate the
breadth of configurations \orcha enables with a case study in which an application configuration is realized on three distinct hardware mappings---a GPU-centric, a CPU/GPU balanced, and a CPU/GPU concurrent
layouts by using different recipes.
\end{abstract}


\maketitle

\section{Introduction}
The end of Dennard scaling \citep{dennard1974design} has led to a paradigm shift
in high-performance computing (HPC) for achieving higher computational capability  by hardware
platforms. Hardware designers are resorting to mechanisms like
specialization that make the platforms harder to program. At the same time,
as computational science tackles ever more challenging problems, the
scientific software is becoming more complex. The programming
models and compiler communities have provided some tools and
mechanisms to alleviate some of the challenges, but what is
increasingly obvious is that there is no magic compiler that can
manage all these challenges.  Different challenges must be addressed
with different abstractions, the abstractions must be able to
interoperate with one another, and tools must enable orchestration of
computation on available resources with cost-benefit analysis for
energy use and time to solution. 

The widely adopted state-of-the-art approaches for tackling heterogeneity are  concerned
with the unification of variants in data structures needed by different
devices. Some of these solutions, in particular those based on 
template meta-programming such as
Kokkos~\citep{CarterEdwardsSunderland14,Trott22}, have had widespread
success. This is because the current generation of machines, including
exascale machines, have most of their compute power in the GPUs. The
other aspect of handling heterogeneity, that of orchestration of data
and computation with task-based
runtimes~\citep{duran2009proposal,heller2013application,legion,bosilca2013parsec,augonnet2009starpu,bosilca2020template},
has not enjoyed as much success. The successful 
ones have been those that have co-developed with their client
application or a class of applications. For existing code bases,
adopting these runtime systems has been very challenging
  primarily because interfacing with these systems requires deep
  refactoring and the code acquires a hard dependency---it
  cannot run without the library becoming available on the target
  platform. Overall, the above techniques adopt a paradigm we call \emph{one-to-many:} meaning many variants of a program result from one abstraction.
  
Our solution for \flashx~\citep{DubeyWeideONeilEtAl22}, 
a Fortran-based multiphysics community code, resulted in \orcha, which reflects direct 
responses to what we perceived as the pain points in handling platform heterogeneity. Unlike the prevalent \emph{one-to-many:} approaches described above, \orcha adopts a \emph{many-to-many} paradigm: each variant of a program results from a specific high-level description of a configuration. What \orcha\ proposes are simple and effective ways of managing the drastic increase in code complexity caused by the need for such different configurations.

We identified five requirements to meet
the needs of the users: (1) a way to permit variants in 
data structures and algorithms demanded by different target devices
without a combinatorial explosion of code versions; (2) a way to
express the map of what to execute where; (3) a mechanism to
orchestrate the movement of data and computation between devices to
execute the map; (4) avoid making the application code acquire a
  hard dependency on the availability of complex tools; and (5) a
reasonable path for gradual adoption of the abstractions without
disrupting the scientific productivity or code 
maintainability. The current implementation of \orcha\ has three stand-alone tools that satisfy the first three requirements. These are: \cgkit\
\citep{RUDI2025107511}, \milhoja\ \citep{ONeilWahibDubeyEtAl22}, and \macro\
\citep{DubeyLeeKlostermanVatai23}, described in detail in our prior work.
Of the three tools, only \macro interacts directly with the maintained source code. It is a short and simple Python script that application developers can maintain themselves. All code required for interfacing with \milhoja 
is generated by a code generator guided by \cgkit.
Therefore, the application code remains entirely ignorant of \milhoja and \cgkit, and it can continue to function without them. This is the key to meeting the fourth requirement.  The code generators are application-specific---either customized or developed from scratch for each application. 
The fifth requirement is met by permitting the application code to selectively use \orcha on some of the code while letting the rest of the code run without it.

The descriptions of \cgkit, \milhoja, and \macro in this paper are included for completeness; the technical rigor of each individual tool has been examined in detail in our prior publications~\citep{RUDI2025107511,ONeilWahibDubeyEtAl22,DubeyLeeKlostermanVatai23}. The present work focuses on the overall design philosophy, the use of code generation as a decoupling layer and a harnessing mechanism among the tools, and demonstration of recipe-driven hardware mapping and the orchestration of example applications.

The contributions of the present work are:
(i)~the overarching design philosophy of \orcha\ and its value to
  the applications
(ii)~the design of code generators for \flashx, the motivating application,
(iii)~the end-to-end realization of the \orcha\ system for \flashx,
and (iv)~demonstration of \orcha\ facilitating multiple configurations of \flashx\ that target different mappings to available hardware.

We demonstrate the use of \orcha\ in \flashx\ using two science
application instances with
different configurations. The first is the Sedov blast wave problem
\citep{sedov1993similarity}, which is one of the standard test cases
for compressible hydrodynamics solvers that can handle strong
shocks. The second is a more complex multiphysics application of cellular
detonation from computational astrophysics that uses a nuclear
network and a computationally expensive equation of state (EOS) in
addition to compressible hydrodynamics \citep{timmes2000}. The second
application also presents an interesting use case where, because of a
third-party dependency that does not have a GPU implementation,
nuclear burning has to be computed on the CPU. And because both
compressible hydrodynamics and nuclear burning need to apply EOS, it
must be computed on both CPU and GPU.  

The paper is organized as follows: \cref{sec:background} provides background,
an overview of \orcha\ is presented in \cref{sec:overview}, the
application-specific interface layer design and  code generation tools are described
in \cref{sec:code-gen}. In \cref{sec:apps}, we describe the use cases
and their different configurations used in the
experiments. \cref{sec:experiments} outlines the experiments and their
results. In \cref{sec:conclusions}, we present our conclusions. 

\section{Background}
\label{sec:background}
Since the introduction of GPUs in HPC, developers have been looking
for ways to utilize them effectively. A key element of being able to
use GPUs is to mitigate the challenge of programming them (e.g.,
\citep{padal2014, Mittal_2015}). The solutions have  taken several
forms that can be broadly categorized into five types. The first were
specialized languages like CUDA and CUDA Fortran, supplied by NVIDIA,
which also provided the first few generations of GPUs for HPC. The
second includes compiler directives-based solutions such as
OpenACC~\citep{openacc} and OpenMP~\citep{openmp}. The third are
domain-specific languages 
\citep{halide,gysi2015stella,clement2018claw,nebo}, several of which
enjoyed brief success in their target communities. The burden of 
keeping up with the growth of software proved to be too great for most
of these communities, and some switched over to the fourth class of
solutions: C++ template meta programming based abstractions such as
Kokkos~\citep{CarterEdwardsSunderland14},
Raja~\citep{BeckingsaleBurmarkHornungEtAl19},
STELLA/GridTools~\citep{gysi2015stella},   OCCA
\citep{MedinaSt-CyWarburton14},   SYCL
\citep{ReindersAshbaughBrodmanEtAl21}, etc.,   where a single
expression of the computation can be specialized to the target device
as needed. These solutions have had the most wide spread success
  and have been adopted by a large number of C++ codes. Their primary
  limitation is that the data structures of the code need to be owned
  by C++, which is not always feasible for legacy codes written in
  Fortran. For Fortran-based codes like \flashx, the only available options are the
directive-based solutions such as OpenMP or OpenACC to compile the code for the GPU. While
necessary, this capability is not sufficient. In a complex code, one
needs the freedom to orchestrate computation and data transfers to maximize
resource utilization.
Additionally, while tools like Kokkos can select between available hardware for executing each function, they do not enable effortless orchestration of high-level control flow with hardware mapping the way CG-Kit, combined with code generation and \milhoja, does.

Some solutions, such as Task-based programming models and/or
tools for shared memory   (e.g., CUDA graphs~\citep{cudagraphs},
OmpSs~\citep{duran2009proposal}) and for distributed memory 
  (e.g., HPX~\citep{heller2013application},   Legion~\citep{legion},
  PaRSEC~\citep{bosilca2013parsec},   StarPU~\citep{augonnet2009starpu},
  TTG~\citep{bosilca2020template}) provide asynchronous data movement
  along with abstractions. These solutions have a critical limitation
  in that they require deep intrusive changes to the code, and
    the code acquires a hard dependency on very 
complex tools that need dedicated, well-funded teams to maintain
and keep up with the changes in platforms.  This becomes a risky
proposition because of a lack of future guarantee. The design
  choice of isolating the tools API from the application and having an
  intermediate layer of code generation in \orcha is meant to
  minimize the negative impact of the lack of future support on the
  applications adopting the tool-chain. As mentioned earlier, the
  applications can be built with or without \cgkit and \milhoja, and
  the extent to which code is converted to macros is at the discretion
  of the developer.
  Importantly, this decoupling does not preclude the use of task-based runtimes themselves: tailored code templates for \cgkit emitting code for one of the portability runtimes mentioned above could occupy the same slot in ORCHA that Milhoja currently fills. The design choice insulates the application from any one runtime, not from task-based runtime solutions.
  \cref{tab:orcha_comparison} compares features and capabilities of \orcha with the other abstractions discussed so far.   

\begin{table*}[t]
    \caption{Comparison of \orcha with existing performance portability approaches. The table highlights differences in abstraction, integration model, the expression of hardware mapping, and orchestration. }
    \label{tab:orcha_comparison}
    \centering
    \scriptsize
    \setlength{\tabcolsep}{3pt}
    \renewcommand{\arraystretch}{1.1}

    \renewcommand\tabularxcolumn[1]{m{#1}}
    \newcolumntype{C}[1]{>{\centering\arraybackslash}m{#1}}
    \newcolumntype{Y}{>{\raggedright\arraybackslash}X}
    
    \begin{tabularx}{\textwidth}{C{2.2cm} Y Y Y Y Y Y Y}
        \toprule
        \textbf{Approach} &
        \textbf{Abstraction} &
        \textbf{Code changes} &
        \textbf{HW Mapping expression} &
        \textbf{Orchestration without source modification} &
        \textbf{Dependency} &
        \textbf{Fit} &
        \textbf{Limits} \\
        \midrule
        
        \textbf{OpenMP / OpenACC} &
        Directive-based parallelism and offload &
        Add directives, manage data regions, minor refactoring &
        In-source directives and clauses &
        Limited (embedded in code) &
        Low–moderate runtime dependence &
        Incremental GPU adoption &
        Limited for global orchestration \\
        
        \addlinespace
        
        \textbf{Kokkos / RAJA} &
        C++ abstractions for execution and data &
        Refactor into policies and data abstractions &
        In-source via policies and layouts &
        Limited (kernel-focused) &
        Strong framework dependence &
        Best for C++ &
        Not suited to other languages, like Fortran \\
        
        \addlinespace
        
        \textbf{Task runtimes (e.g., Legion, StarPU)} &
        Task graphs with runtime scheduling &
        Restructure into tasks and data regions &
        Runtime APIs and task graphs &
        Yes (via runtime model) &
        Strong runtime dependence &
        Flexible scheduling &
        High integration cost for legacy codes \\
        
        \addlinespace

        \rowcolor{gray!10}
        \textbf{ORCHA} &
        Recipes + code generation + lightweight runtime &
        Define APIs, optional macro usage, minimal data annotations &
        External recipes + generated code &
        Yes (recipe-driven) &
        No dependency in maintained source; isolated in generated code &
        Enables variant exploration and mixed HW mappings while keeping the source code unchanged &
        Requires application-specific code generator development \\
        
        \bottomrule
    \end{tabularx}
\end{table*}

  The final set consists of  HPC
  languages such as Chapel~\citep{chapel} or Co-array
  Fortran~\citep{coarray}. These languages are well designed but suffer
  from a lack of  enough resources to make them viable for
    production codes. Chapel is not supported by vendors, and while 
  many features of Co-array Fortran
have been incorporated into the standard, its parallel model, PGAS \citep{almasi2011pgas}, is incompatible with MPI \citep{mpi40}, where vendors put most of their communication optimization efforts.

\section{Overview of \orcha}
\label{sec:overview}
Contemporary computational science is faced with two axes of
simultaneously increasing heterogeneity. While hardware heterogeneity is
well known and understood, application software heterogeneity is not
yet on the radar of the HPC community. Science models tend to
become more complex with time, introducing significant variations in
data structure and data movement requirements within the
application. Additionally,  all components of an application may not
be suitable for all variants of hardware architecture in a
platform. To effectively utilize available resources, applications
need: (1) data structures and algorithms suitable for target devices,
(2) a way to conceptualize and describe a map of computation to target
devices, and (3) a way to execute the map by moving data and
computation to devices efficiently. Given the complexity and size of
the cutting-edge scientific applications, the above requirements must
be met without having to redesign algorithms or refactoring existing
code for each new generation of platforms. Maintaining too many
different variants of the implementations is also not a viable option
from the maintainability and productivity standpoint. Therefore, the
abstractions must be designed so that the application can be tuned for the
target platform with minimal changes to the source code. As mentioned
earlier, for C++ codes, Kokkos has been very successful in this regard for the first our of
three requirements discussed above. One of the tools in \orcha, a
customized \macro, addresses the first concern for non-C++
codes. The other two tools in \orcha\ are designed to address the
features missing in the current widely adopted solutions -- the ability to
unify algorithmic variants and adopt different orchestration of data
and computation for different architectures without changing the
source code. Here, we briefly describe the main tools for
completeness; for more details, readers are referred to
\citep{RUDI2025107511,ONeilWahibDubeyEtAl22,DubeyLeeKlostermanVatai23}. 

It is worth noting that the three tools described in the following sections (\cgkit, \milhoja, and \macro) are the reference toolchain of \orcha. Each of these tools can be replaced by external tools of similar capability, by modifying the code generator accordingly. The overarching design philosophy of \orcha, a recipe-driven hardware orchestration based on code generation, remains applicable in principle, as long as the overall design is carried out.

\subsection{\cgkit}
\label{sec:cgkit}
\cgkit\ enables knowledgeable users to express execution control flow and to map computations to specific hardware targets without modifying the core Fortran or C++ source code of the application. The control flow and computation are expressed in recipes, which describe how code variants are composed and executed.

Our goal is to give application developers a practical way to express implementation variants and explore diverse optimization strategies. A problem arises, for example, when porting an application to a new platform. Developers must identify performance-critical sections and optimize them for the new target architecture. This often results in multiple algorithmic variants, which in turn leads to extensive branching, making the code difficult to maintain and reason about. 
 
\cgkit\ addresses this challenge by implementing the variants as labeled templates containing code snippets with specific implementation choices. Recipes then compose these templates into complete source files, describing algorithmic variants at a high level. 
To generate source code, CG-Kit translates the recipe into a statically optimized directed acyclic graph (DAG) and constructs a tree representation of the code to be emitted. While abstract syntax trees (ASTs) could be used for the latter, they are not well suited for human inspection. A key design goal of CG-Kit is that all generated code and intermediate representations remain human-readable, enabling developers to better understand, debug, and reason about performance. To this end, CG-Kit uses parameterized source trees (PSTs), a human-friendly alternative to ASTs inspired by template-based code generation systems.
A PST is a tree where leaves are parametrized code snippets written in the target language (e.g., C++ or Fortran) rather than language tokens, so the intermediate representation and the final source tree can be read as ordinary source code. Full details of the PST representation and its composition rules are given in~\citep{RUDI2025107511}.
In this approach, the templates described earlier become building blocks of PSTs. 
The final PST is constructed by composing these template PSTs as dictated by the optimized recipe graph which is then parsed to emit compilable source code in the original programming language.

CG-Kit is implemented in Python and operates as a text-based source transformation tool. As a result, it is language-agnostic and can be applied to source code written in any programming language.

\subsection{\milhoja}
\label{sec:milhoja}

\milhoja\ is a domain-specific runtime library written in C++,
designed to execute a fixed execution graph of computational operations repeatedly. 
It is useful in applications whose control flow on a compute node
remains unchanged from one evolution step to the next (e.g., time stepping in a physics simulation) throughout its
execution.  By design---instead of being feature rich---\milhoja is implemented with minimum possible capabilities to asynchronously execute computations and transfer data to/from device as dictated by the graph emitted by \cgkit. Its core communication engine is an extended finite-state machine designed to move packets of data placed in a queue to their designated destination and execute computations using thread teams. 

For a given simulation, \milhoja operates with thread teams that are created when the simulation starts and persist until
the simulation terminates. The thread teams are run in cycles such that for each cycle, a team is assigned a dataset and work to be done
on that dataset. The pairing of the work to be done with the dataset is viewed as a task, so that thread teams implement task-based parallelism within a node via  thread-level parallelism. 
Note that the tasks are strictly devoid of any distributed memory parallel constructs. Thus, a task never has an MPI operation embedded in it. If a solver requires MPI operations in the middle of their computations, it must be broken up into $n+1$ operations where $n$ is the number of
invocations of MPI calls in the solver.

The current implementation of \milhoja is accompanied by a code generator customized for block-structured data. Applications that do not use block-structured data will need a code generator tailored for their specific data model.

\subsection{\macro}
\label{sec:macro}
Numerical algorithms typically have arithmetic interspersed with
control flow, and their arithmetic is expressed in terms of their data structures. If we abstract away the knowledge of data structure when expressing the arithmetic, the source code results in fewer variants of the code to be maintained. In C++-based
solutions, this is done through templating. In \orcha, we make it work by associating a unique string with each code-snippet that we wish to use as a building block of code, and/or customize for different device architectures. We call these code-snippets macros because they are similar to C preprocessor macros in how they are defined. They permit arguments, are allowed to be inlined in
a regular programming language statement, and one macro name can be
embedded in another macro’s definition as long as self-reference is
avoided. Our \macro differs from the preprocessor in that it permits multiple alternative definitions of macros and includes the logic to arbitrate on selection of a specific definition. This feature enables a great deal of flexibility such as redefining partial sets of definitions and allowing inheritance branching.

This mechanism
permits the unification of code variants where the differences are in
data layout and movement rather than the arithmetic of the
calculation. Macros also provide a way to support different
loop-level parallelization options within the static code by having
those constructs be alternative definitions for the macros that can be
inserted into the source code instead of the pragmas or other
constructs. Macros also allow the code
to be decomposed into components of arbitrary granularity, facilitating
effective orchestration and reuse of code snippets.

\section{Code Generators}
\label{sec:code-gen}
The tools described in the previous section are the core stand-alone tools
designed for specific tasks that many applications need. They expect
input to be provided to them in a specified format, and they perform
the required work based on that input.
In this regard, they are similar to several other tools described in
\cref{sec:background}, that can, in principle, replace tools in \orcha.
To make the end-to-end solution work, applications need the tools
to be stitched together into a configuration and code generation pipeline 
so that they can interoperate with one another in an application-specific way.

We introduce a few definitions to facilitate the discussion of how \orcha operates.
\begin{itemize}
    \item A \taskfunction is a collection of an application's API functions that are executed as a group on a target device.
        For example, in \flashx, a user may opt to have hydrodynamics and gravity computed on the GPU,
        while nuclear burn is computed concurrently on the CPU.
        In this instance, hydrodynamics and gravity together form one \taskfunction, because they are computed on the same device without any data transfer in between.
        A second \taskfunction will be formed by nuclear burn alone because it is to be executed on a different device.
    \item A \dataitem is the container that supplies one \taskfunction with all the data it needs to execute. It is generated in one of two shapes depending on where the \taskfunction runs. When the \taskfunction executes on a device, the \dataitem is generated as a \datapacket---a flattened structure that bundles the required data into a single transfer to the device. On the other hand, when the \taskfunction executes on the host, the \dataitem is generated as a \tilewrapper---a lightweight wrapper that exposes host-resident data to the \milhoja runtime. The distinction reflects the shape of the generated code, not an asymmetry between host and device data characteristics. In the above example, the data for hydrodynamics and gravity (executed on GPU) is bundled into a \datapacket-shaped \dataitem, while the data for nuclear burn (executed on CPU) is wrapped in a \tilewrapper-shaped \dataitem. From hereon, we use \dataitem generically and refer to the specific shapes only for distinguishing the two code-generation paths.
    \item Static Code is the source code that implements the API functions of the application.
        This code has embedded application-specific macros that are expanded by the \macro.
\end{itemize}

We need two major classes of application-specific generated code for
using \orcha: one that queries the application and consolidates all
the data needed to write a \dataitem, and another that combines all 
the high-level API functions into a \taskfunction\ as dictated by the 
recipe and translated by \cgkit. Since the high-level API functions may rely on application data, there is an inherent co-dependency
between a \taskfunction and a \dataitem.
So, once a \taskfunction is determined, the corresponding \dataitem
 is also fully determined.

To generate application-specific code, the tools need to (1) know the hardware requirements for 
running a simulation, (2) translate the existing static code into a format that
is understood by the code generation tools, (3) generate code that
calls each function or subroutine at the top of the call stack in the
\taskfunction, as well as code that can transfer data between devices
as needed. For \flashx, we have opted to use annotations in the API of
the included code units as a way of communicating
static code requirements to
\orcha\ code generation tools. We impose the following
  requirements on the static code:
\begin{itemize}
  \item Static code is thread safe.
  \item Any function to be accessible to a recipe must be promoted to
    the API of the static code.
  \item Static code does not allocate any scratch space; it receives
    needed scratch space as an argument in its interface.
  \item Only constants are allowed to be global; all other variables
    must be arguments in the interface.
  \item All non-scalar scratch space, along with its size, is described in the
    annotations.
  \item All the data that is needed for the computation is enumerated
    in the annotations.
\end{itemize}
Except for the annotations, all other requirements are consistent with
any programming model that one might use for moving data and
computation between devices; therefore, \orcha does not place any unusual 
burden on the applications' static code. More importantly, these modifications do not tie the
application to any specific tool.
Those assumptions highlight the major difference in \orcha's approach relative to existing solutions. With \orcha, an existing application code is \emph{not} refactored to fit a particular portability framework. Instead, the application needs to have clear interfaces as defined above---adopting universal best practices. The goal of \orcha\ is to utilize these interfaces and to employ code generation as well as a simple runtime library to adapt the application to different platforms.

We define a contrived example to illustrate how the tools of \orcha
and the code generators work together to create an instance of
an application from a user-defined recipe. We assume four subroutines,
\texttt{P1}, \texttt{P2}, \texttt{P3}, and \texttt{P4}, in the code API that are accessible to the recipes.
\cref{lst:contrived_defs} shows the signatures of four subroutines and the implementation of \texttt{P1}.
The code annotations, starting with double exclamation marks, \texttt{!!},
explain each subroutine argument to the code generators,
thereby determining the required data for subroutines.
In practice, the subroutine annotations are structured in a more confined,
systematic form to be easily parsed by the code generators, but we use plain text in this example for explanatory purposes. 
The implementation of \texttt{P1} in \cref{lst:contrived_defs} includes macros,
which are marked with \texttt{@M} before their names.
The reader can assume similar implementations for \texttt{P2}, \texttt{P3}, and \texttt{P4}.

\begin{lstlisting}[
    language=Fortran,
    caption={Example code APIs to be used in a recipe for \orcha. There are four subroutines, \texttt{P1}, \texttt{P2}, \texttt{P3}, and \texttt{P4}. The argument annotation block and subroutine implementation with macros are expanded only for \texttt{P1}, for simplicity.},
    label={lst:contrived_defs}
]
!! Ustate: input/output data structure described by Udesc
!! Scr1, Scr2: scratch spaces described by Udesc
!! S1, S2, S3: input scalars
!! Udesc: it includes the size and dimensionality specification for the data used in the example
subroutine P1(Ustate, Udesc, Scr1, Scr2, S1, S2, S3)
    @M Udesc_declare
    real, dimension(@M Udesc_view), intent(inout) :: Ustate, Scr1, Scr2
    real, intent(in) :: S1, S2, S3
    integer :: i, j, k
    @M loop_ijk_begin
       @M Comp1(Ustate, Scr1, Scr2, S1)
    @M loop_ijk_spl_end
    @M loop_ijk_spl_begin
       @M Comp2(Ustate, Scr1, S2, S3)
    @M loop_ijk_end
end subroutine P1

subroutine P2(Ustate, Udesc, Scr, S1, S2)
subroutine P3(Ustate, Udesc, Scr, S1)
subroutine P4(Ustate, Udesc, S1, S2, S3)
\end{lstlisting}

In \cref{lst:contrived_defs}, we assume that all macros except \texttt{Udesc\_declare} and \texttt{Udesc\_view}
have two sets of definitions: one for CPU and one for GPU. The macros
\texttt{loop\_ijk\_spl\_begin} and \texttt{loop\_ijk\_spl\_end} have null definitions for
the CPU, while for the GPU, they are
responsible for breaking the loops,
so that on the CPU, both \texttt{Comp1} and \texttt{Comp2} are
in one kernel, whereas on the GPU, they are on separate kernels.
In this way, \macro expands \texttt{P1} to two different variants for two target hardware,  \texttt{P1\_CPU} and \texttt{P1\_GPU}.
The actual example macro definitions and the final expanded versions of \texttt{P1} are provided in \cref{appendix:macro}.

We can now write several different recipes, each of which orchestrates the
computations differently. For example, we could compute \texttt{P1} on
the GPU, followed by \texttt{P2} and \texttt{P3} computed on the GPU,
while \texttt{P4} is computed on the CPU simultaneously.
The invocation order, dependencies, and hardware mapping of subroutines
are described in recipes; thus, several other application configurations
are possible by changing recipes, without redefining the application's APIs.

\begin{lstlisting}[
    language=Python,
    caption={Example recipe using \recipetools, utilizing APIs defined in \cref{lst:contrived_defs}.},
    label={lst:example_recipe}
]
import FlashX_RecipeTools as flashx

# Initialize recipe
recipe = flashx.TimeStepRecipe()
begin = recipe.begin_orchestration(after=recipe.root)

# Add nodes
p1 = recipe.add_work("P1", after=begin, map_to="CPU,GPU")
p2 = recipe.add_work("P2", after=p1, map_to="GPU")
p3 = recipe.add_work("P3", after=p2, map_to="GPU")
p4 = recipe.add_work("P4", after=p1, map_to="CPU")

# End recipe
end = recipe.end_orchestration(begin_node=begin, after=[p3, p4])
\end{lstlisting}

\begin{figure}
    \centering
    \includegraphics[width=0.8\linewidth]{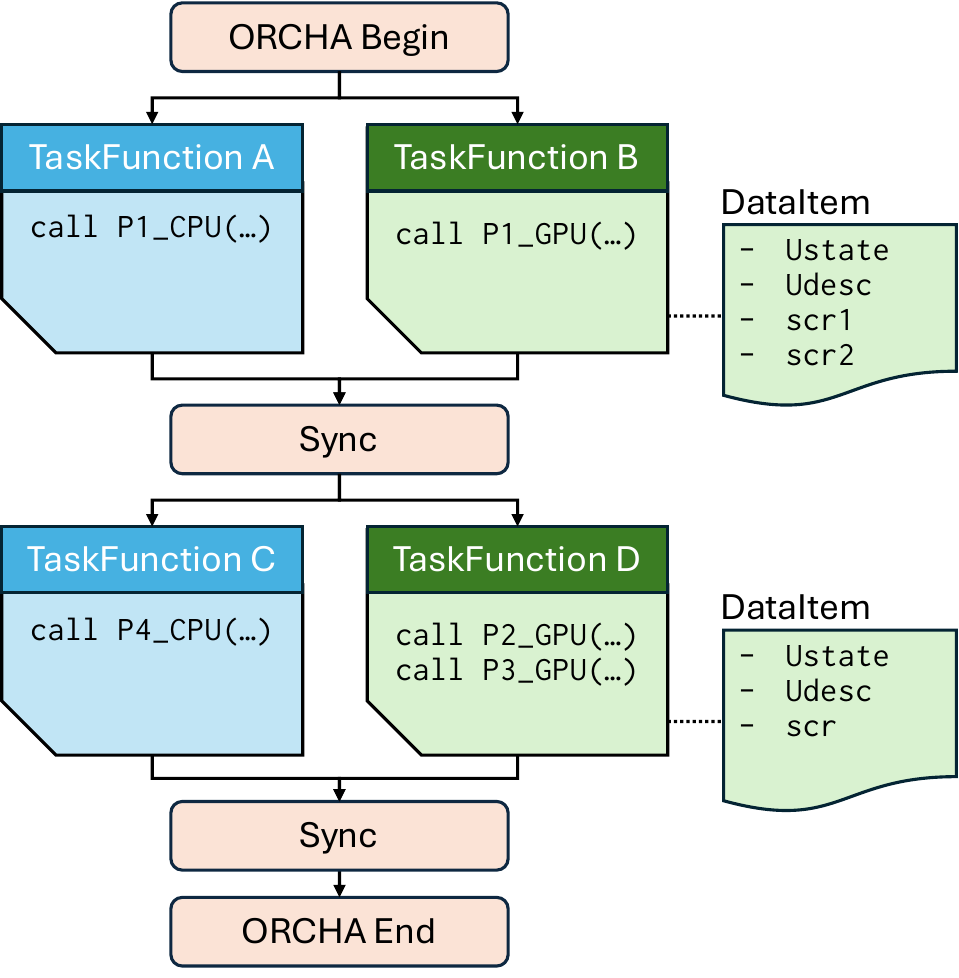}
    \caption{The \taskfunction flow graph corresponding to the example recipe, \cref{lst:example_recipe}, parsed and compiled by \recipetools.
    Green boxes represent the GPU \taskfunction{}s,
    while blue boxes represent the CPU \taskfunction{}s.
    The \dataitem{}s for the GPU \taskfunction{}s (hence, \datapacket{}s) are also depicted,
    which contain unions of data that need to be transferred
    to the GPU to call the subroutines contained in the GPU \taskfunction{}s.
    The \dataitem{}s for the CPU \taskfunction{}s are omitted in this figure for simplicity.}\label{fig:tf_flow_graph}
\end{figure}

\cref{lst:example_recipe} presents an example recipe
that uses the subroutines defined in \cref{lst:contrived_defs},
assuming they are \flashx’s physics APIs.
The recipe is constructed using a Python tool called \recipetools,
which provides a class for building recipes and code generators
for the target application, \flashx.
The code of \recipetools\ utilizes \cgkit; it simplifies recipe writing for \flashx.
For the recipe in \cref{lst:example_recipe}, \texttt{P1}
will be computed on both the CPU and the GPU, followed by computing \texttt{P2} and the  \texttt{P3} on the GPU,
and simultaneously computing \texttt{P4} on the CPU. The corresponding 
flow graph 
is shown in \cref{fig:tf_flow_graph}.

In the figure, because \texttt{P1}  needs to be executed on both the CPU and the GPU,  two distinct names \texttt{P1\_GPU}, and \texttt{P1\_CPU} are needed for different compiled object files. These object files are to be assigned to two created \taskfunction{}s A and B respectively. A \tilewrapper, a \dataitem for CPU \taskfunction, needs to be created for \taskfunction A and
a \datapacket, a \dataitem for GPU \taskfunction, is needed for \taskfunction B. \taskfunction{}s C along with a \tilewrapper will be created for computing \texttt{P4} on the CPU , and \taskfunction D along with a \datapacket will be created for computing \texttt{P2} and \texttt{P3} on the GPU. \texttt{P2} and \texttt{P3} are put into the same \taskfunction and their data into the same \datapacket because no data transfer is needed between their execution. Additionally, \taskfunction{}s C and D can be executed in parallel based on the dependency analysis done by \cgkit.

\begin{figure*}[hbt!]
    \centering
    \includegraphics[width=\textwidth]{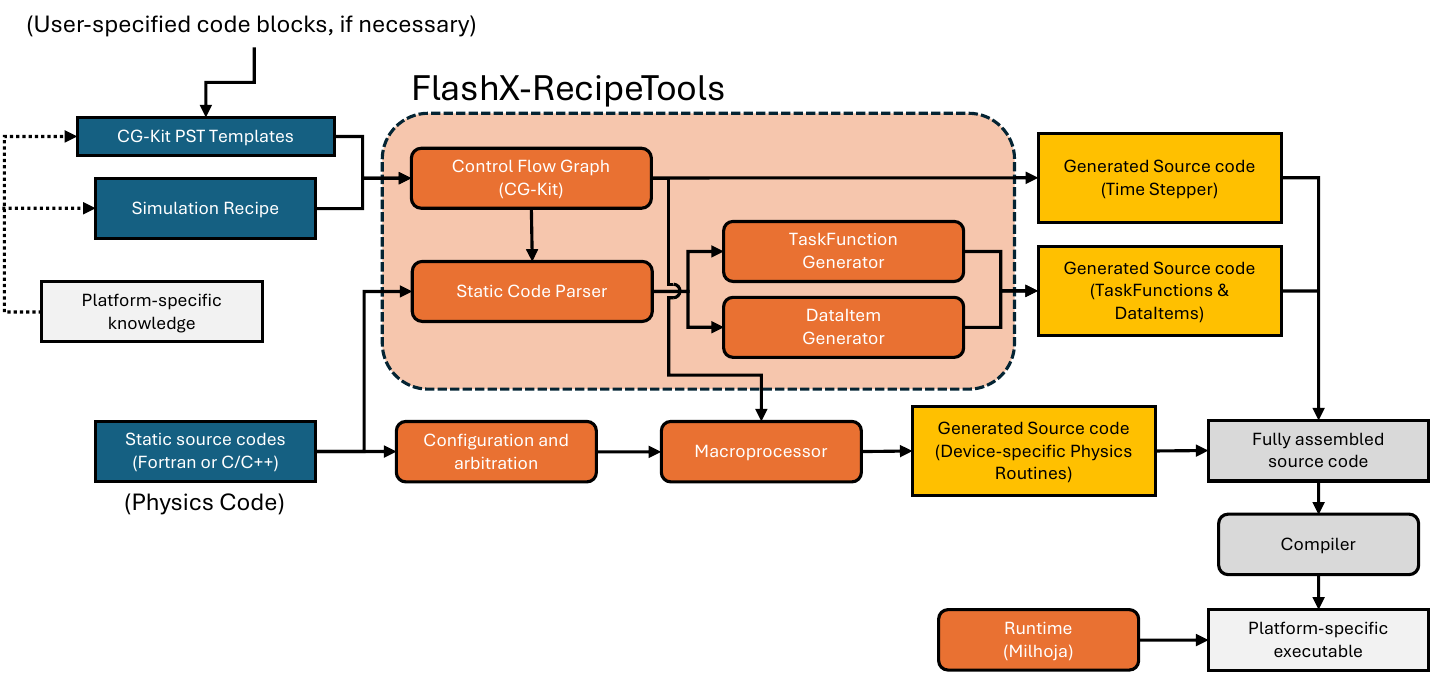}
    \caption{Overview of the code generation pipeline for \flashx with the new \recipetools (highlighted in the center by the large orange box) and \orcha.
    \recipetools parses the user-written simulation recipe,
    determines the control flow graph, and then generates the source codes
    for the platform-specific target application.
    Internally, it extracts information from the given recipe
    and interacts with three different tools: \cgkit~\citep{RUDI2025107511}, \milhoja~\citep{ONeilWahibDubeyEtAl22},
    and \macro~\citep{DubeyLeeKlostermanVatai23}.}
    \label{fig:pipeline}
\end{figure*}

To generate the \dataitem{}s, static code parser parses the annotations in all relevant APIs and consolidates data for  each \taskfunction. 
The \dataitem{}s contain only non-scalar values, assuming input scalars are trivially initialized on the target device.
The \datapacket for the \taskfunction D contains a union of data required by
subroutines \texttt{P2} and \texttt{P3}, allowing one scratch space to be
recycled by both subroutines. The code for creating \taskfunction{}s and \dataitem{}s is generated by the \recipetools as illustrated by the pipeline is shown in \cref{fig:pipeline}. We next describe code generation done by the \recipetools.

\subsection{Static Code Parser}
The static code parser limits its parsing to the annotations in the
code, therefore, it is applicable to source code in any language as long as the syntax of the annotations is the same. Therefore, we will assume the static code to be in Fortran without any loss of generality.

The code generators in \orcha need a complete inventory of all data that needs to be associated with each \taskfunction and the corresponding \dataitem. %\datapacket{}/\tilewrapper.
In \flashx, each physics unit publishes its API, which includes all the interfaces that are accessible to recipes and other units in the code. 
Recall that to be usable in a recipe the interfaces must include all scratch space to be used as a part of the argument list,
therefore, all required data except global constants can be captured by enumerating the APIs' input arguments.
We use a structured syntax for the annotations which also includes meta-information about each of the arguments. For instance, this information may indicate whether an argument is a state variable of an AMR block, scratch space with a specific size and starting indices, or a simple scalar. 
 
The annotations are parsed by 
the \recipetools and converted to the JSON format.
The JSON files are consumed to create a \taskfunction Specification,
which is a representation of unified collection of subroutines and their data. The \taskfunction Specifications are in turn 
 passed to the \taskfunction and \dataitem generators
to generate the source code.

\subsection{\taskfunction\ Generator}
The \taskfunction\ generator uses the \taskfunction
Specification to generate code for (1) an interoperability
layer between C++ and Fortran, (2) extracting the data from the \dataitem, and (3) wrapper code that calls
each subroutine included in the \taskfunction. %Specification from the recipe tool. 
The C++ to Fortran interoperability layer is
needed because the vast majority of static code in \flashx\ is written in
Fortran, whereas \milhoja\ is written in C++. 
The generated interoperability layers are responsible for extracting the data in two steps.
In the first step, a C function
extracts data from the \dataitem\ into  C pointers or primitive
values. Then, a Fortran subroutine converts each C pointer into a
Fortran pointer or a Fortran primitive value. Once
this conversion is complete, the data context is effectively ready for
the \taskfunction. 

\subsection{\dataitem\ Generators}\label{subsec:datapacket_generator}
The \datapacket\ Generator is one of two \dataitem generators. It is responsible for writing code that combines any data described in a \taskfunction specification  to be executed on a target device
into a flattened data structure that can be moved as a single data packet. The data extracted from \datapacket{s} at the
device end sets up the structure of the data in the device memory,
requesting the necessary memory from \milhoja, and moving pointers
into pinned memory and back.

\datapacket Generator can also specialize application-specific data structures
for better efficiency.
For example, an AMR-based application like \flashx prefers to concatenate
several AMR blocks into one \datapacket to 
minimize the data transfer latency since individual blocks are too small for effective utilization of GPUs.
\recipetools 
 can expand any block-specific data to include $n$ AMR blocks, where $n$ is provided by the users at runtime to adjust for their specific use case.  
Similar \datapacket specialization can be achieved for other applications
by providing different code templates to the \datapacket Generator. \cref{sec:experiments} demonstrates its use in \flashx.

The \tilewrapper Generator is the other \dataitem generator. It 
works similarly to the \datapacket Generator, creating a simpler \dataitem called a \tilewrapper. \taskfunction{s} to be performed on the host use this instead of \datapacket{s} because there is no data movement between devices. Therefore, the \tilewrapper Generator only needs to wrap data specified in the \taskfunction Specification so that \milhoja can understand it, instead of combining chunks of data into a flattened data structure.

\subsection{Code Generator Development Effort and Generalization}

\begin{table*}[t]
    \caption{
        Effort required to adapt \orcha to example applications with different discretizations.
        Other application classes follow similar trends, with effort determined primarily by
        data layout and access patterns differences.
    }
    \label{tab:effort}
    \centering
    \scriptsize
    \setlength{\tabcolsep}{3pt}
    \renewcommand{\arraystretch}{1.1}

    \renewcommand\tabularxcolumn[1]{m{#1}}
    \newcolumntype{C}[1]{>{\centering\arraybackslash}m{#1}}
    \newcolumntype{Y}{>{\raggedright\arraybackslash}X}
    
    \begin{tabularx}{\textwidth}{C{2.2cm} Y Y Y Y Y Y}
        \toprule
        \textbf{Adaptation Efforts} &
        \textbf{Uniform Grid} &
        \textbf{Structured AMR} &
        \textbf{Unstructured AMR} &
        \textbf{Dense Linear Algebra} &
        \textbf{Sparse Linear Algebra} \\
        \midrule

        \textbf{\milhoja runtime} &
        Reusable &
        Reusable &
        Requires extension for data packing &
        Requires extension for data packing &
        Requires extension for data packing \\
        \addlinespace

        \textbf{Static Code Parser} &
        Mostly reusable &
        Mostly reusable &
        Requires extension for data access patterns &
        Requires extension for data access patterns &
        Requires extension for data access patterns \\
        \addlinespace
        
        \textbf{\taskfunction and \dataitem generators} &
        Reusable &
        Reusable &
        Mostly reusable, but may require a new data model &
        Mostly reusable, but may require a new data model &
        Significant redesigning may be needed \\

        \bottomrule
    \end{tabularx}
\end{table*}

\cref{tab:effort} summarizes the minimum effort required to adapt \orcha code generators across different classes of applications. Currently available code generators in \orcha
assume a block-structured data model. Apart from this assumption, the \taskfunction and the \dataitem generators are application-agnostic. Static code parser is the only application-specific component of the toolchain because it has to know about the data shape, access patterns, and interface metadata of the application. As a result, applications that either use or can mimic the block-structured data model can moslty reuse the \taskfunction and the \dataitem generators. A significant portion of the static code parser can also be reused by such applications with changes relating to how the meta-information about the API of the application is communicated to the rest of the toolchain.

In contrast, applications with different data model are likely to require non-trivial amount of work 
in data packing/unpacking logic in \milhoja. Note that 
the required effort depends upon the degree of alignment of the application's data model and iteration structure rather than the size or the number of computational kernels in target applications. This suggests that \orcha adoption efforts can be amortized across families of applications
sharing similar data layouts and execution patterns. 

\section{Applications}
\label{sec:apps}

We demonstrate \orcha\ with two application configurations from  \flashx   that are included in the distributed code as example use cases.  
To help in understanding the applications and the experiments we give a brief overview of \flashx first.

\subsection{\flashx}
\label{sec:flashx}
The \flashx\ code is a component-based software system for building
multiphysics simulation applications that are formulated largely as a
collection of partial and ordinary differential equations. As with any
scientific application addressing complex multiphysics phenomena,
there is also support for other types of mathematical models, including
algebraic equations, table lookups, and parameterized models. In
\flashx, different permutations and combinations of various components
generate different applications for multiple science
domains. Components are allowed to have multiple alternative
implementations to account for either different physical regimes, or
different fidelity requirements, in an application
instance. The components are self-describing in that they include
meta-information about how they will be included in various
application configurations. A tool {\em setup} parses the
meta-information, arbitrates on which components and/or which of a
component's implementations to include, and assembles an application instance. 

The primary discretization in \flashx\ is finite-volume,
though some solvers use finite-difference. The discretized mesh uses
structured adaptive mesh refinement (AMR) to enable higher resolution
only where needed as a way to optimize required resources for memory
and computation. The physical domain is decomposed into {\em blocks}
where each block has the same number of cells and all the cells in one
block have the same physical size. The physical size of the cells in a
block depends upon its level of resolution, so blocks at different
levels span different amounts of physical domain. It is possible to
further decompose a block into smaller chunks, {\em tiles}, if
necessary.  

The mesh is Eulerian, but an overlapping Lagrangian
framework also exists that can exchange physics quantities with the
mesh. Physics solvers available with the current distribution include
compressible and incompressible hydrodynamics, gravity in various
flavors, several equations of state, source terms for nuclear 
burning and chemical transport, support for multiphase flows, and
fluid-structure interaction. Some capabilities, such as neutrino
transport \citep{thavappiragasam2024performance}, 
can be imported as modules and interoperate
with the code. 

\subsection{Hydrodynamics Only Application - Sedov Blast Wave}
\label{sec:sedov}
The Sedov blast problem is a pure hydrodynamics (Hydro) application where the
domain is initialized with a pressure spike at the center. That spike
results in a shock wave that moves out in a radial direction, moving
away from the center. It is one of the famous test cases for compressible
hydrodynamics solvers because an analytical expression exists for how
far the shock has traveled in a given time. This application acts as a
baseline case for our experiments because the only other physics
needed for this application is the ideal gas gamma law EOS, which is
a simple algebraic expression that adds negligible work to the
computation. Therefore, for all practical purposes, it is considered a single physics
application. 

\subsection{Multiphysics Cellular Problem}
\label{sec:cellular}
The Cellular Nuclear Burning problem \citep{timmes2000} uses
compressible hydrodynamics, an equation of state for ionized
plasma formulation using Helmholtz free energy (H-EOS), and a 13
isotope $\alpha$-chain plus heavy-ion reaction network for nuclear
burning (Burn). Thus, there are three physics solvers, each of which
has non-negligible contribution to the computations. The two-dimensional
calculations are  performed in a planar geometry of size $256.0$\,cm by
$25.0$\,cm.  The initial conditions consist of a constant density of
$10^7$\,g\,cm$^{-3}$, temperature of $2 \times 10^{8}$\,K, composition of
pure carbon $\text{X}(^{12}\text{C})=1$, and material velocity of $v_{x}=v_{y}=0$\,cm\,s$^{-1}$.  Near the $x=0$ boundary, the initial conditions are
perturbed to a density of $4.236 \times 10^7$\,g\,cm$^{-3}$, temperature
of 4.423$\times$10$^9$\,K, and material velocity of $v_{x} =
2.876 \times 10^8$\,cm\,s$^{-1}$. The initial conditions and
perturbation given above ignite the nuclear fuel, accelerate the
material, and produce an over-driven detonation that propagates along
the $x$-axis.  After some time, depending on the spatial resolution
and boundary conditions, longitudinal instabilities in the density
cause the planar detonation to evolve into a complex, time-dependent
structure.  

We selected this application because it includes three computationally heavy physics solvers that impose different constraints on how work can be divided among resources on a platform. It can be configured in many different ways to highlight how end users can explore
different options to orchestrate their computations, including any cost-benefit trade-offs. Additionally, the Burn network has dependency on a third party library that is not GPU-ready. Therefore, its 
inclusion highlights an
essential feature of \orcha; the ramp-on path for gradual adaptation of code components for heterogeneity.  Hydro can be run in one of three modes; CPU-only, GPU-only, or by splitting the work between the two. Both Hydro and Burn need to apply H-EOS, which
implies that if any part of Hydro is running on the GPU, we need to
compile H-EOS for both CPU and GPU. Another degree of freedom that we
can explore with this application is that, at the cost of some loss in
solution fidelity, it is possible to have the computations of Hydro
and Burn proceed in parallel on different devices. The loss in
fidelity arises because Burn updates the abundances of nuclear species
based on nuclear reactions, while hydro updates them through
advection. Here, we are using 2-stage Runge-Kutta time integration in
Hydro, where the contribution from Burn at time step $t-1$ is advected
in the first stage at time $t$ in Hydro, but any advection done in the
second stage is obliterated when contribution from Burn at time step
$t$ is added to the solution. Because the contribution from Burn is
typically more significant than through advection, the error introduced is
usually small. When it comes to making trade-offs between time to
solution and fidelity, such considerations can become important.

Because the current generation of
platforms that can be used to run experiments are still CPU/GPU
dominant, we have only used orchestration between CPU and
GPU. However, 
that does not imply any limitation in \orcha's  ability to utilize other 
accelerators in the future as they become available.

\section{Experiments}
\label{sec:experiments}
The experiments described in this section were performed on
Perlmutter at NERSC. Perlmutter is a Hewlett-Packard Enterprise Cray
EX heterogeneous supercomputer with 3,072 CPU-only and 1,792
GPU-accelerated nodes. The CPUs are AMD EPYC with 64 cores per
CPU. The CPU nodes have two CPUs per node, while each accelerated node
has one CPU and 4 NVIDIA A100 GPUs. All of our experiments were
conducted on one accelerated node of the machine since node-level
heterogeneity is what \orcha\ addresses. We run the experiments in
several modes. The easiest mode is to run only on the CPU, utilizing
all the cores available on the node. We use this as the reference
mode. 

The Hydro computation can have either $1$ or $2$ MPI operations
interspersed with local node computations, depending on the choice of
high-level control flow.  The Hydro solver, \spark\ \citep{Couch2021}, is an
explicit finite-volume method where the guardcells form a halo around
the block of cells to be updated. The time integration uses a 2nd-order Runge-Kutta method where the guardcells need to be refreshed for
each stage. We assume that the first refresh occurs before we invoke
\orcha. The second refresh can be avoided if we use a communication
avoidance mechanism where we make the halo twice as thick as it needs
to be and redundantly compute the inner portion of the halo in the
first stage as though it was part of the domain section to be
updated. This trick avoids MPI communications needed to fill the halo
of guadcells between the first and second stages. Additionally, in AMR, flux-correction is needed when two adjacent blocks are at different refinement levels to reconcile physical quantities at their boundaries. Since these blocks may belong to different MPI ranks, a global communication step is necessary. \milhoja only handles computations that do not involve any MPI operations, therefore, computation of Hydro cannot be completed in a single task
function. One can do it with two \taskfunction{s}, one that does
everything, including updating the solution using communication
avoidance, and the second that updates the values corresponding to the
cells that are needed for flux-correction. However, this second
operation is too low in computational intensity to be worth sending to
the GPU, and we typically complete it outside of \orcha.

We should note that this section is devoted to demonstrating \orcha's
core value---flexible hardware exploration 
rather than showcasing peak performance. 
In the two representative multiphysics tests, we generated multiple distinct builds
with different hardware configurations by simply tweaking the simulation recipe.
Although the current prototype lacks kernel- and runtime-level optimization,
we still compare \orcha-generated builds against a hand-tuned CPU baseline.
The underlying physics routines 
are essentially the same that 
have been employed successfully on CPU-only machines for several years. They differ only in including OpenACC compiler directives expanded by \macro to generate GPU kernels. We expect more significant performance improvement once the physics routines are refactored to become more GPU-friendly, e.g., improving cache locality, minimizing register usage, and eliminating automatic arrays in the internal physics routines.
Irrespective, these preliminary results demonstrate how \orcha enables users to  effortlessly regenerate and evaluate alternative hardware mappings,
allowing them to determine the optimal configuration for their production.

\subsection{Sedov Blast Wave Results}
For the Sedov example, the options for offloading are (1) all of
Hydro+EOS computations to the 4 GPUs using 4 MPI ranks, and each rank
talking to one GPU, or (2) operating in a hybrid mode where two thread
teams are in operation with different portions of data assigned to the two thread teams.
One thread team runs on the CPU with multiple threads while the other handles the execution on the GPU.
The number of thread teams
and the work-split ratio can be controlled at runtime.

\begin{figure}
    \centering
    \includegraphics[width=0.9\linewidth]{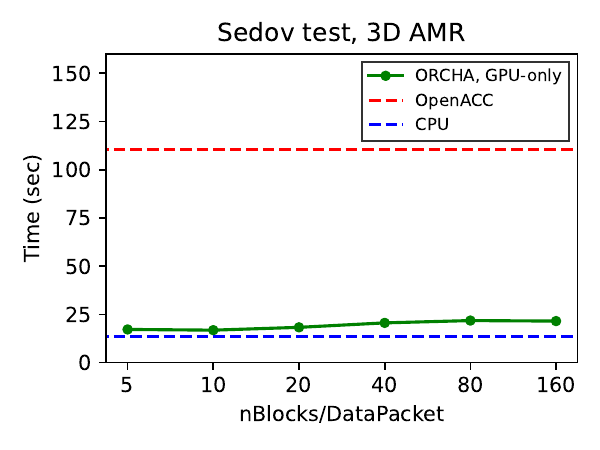}
    \caption{Performance comparison for the Sedov blast application. The two dotted lines are a single data point and are used as references. The red dotted line shows the performance when OpenACC is used without using \orcha, while the blue dotted line shows the performance on CPU with 64 MPI ranks (1 rank per core). The solid green line shows computational time on the GPUs using \orcha in GPU-only mode, as the number of blocks consolidated into a single \datapacket\ is increased from 5 to 160 blocks ($x$-axis in log-scale). 
    }\label{fig:sedov3d-perf}
\end{figure}

\cref{fig:sedov3d-perf}
shows the performance using \orcha\ (green solid line) compared against two reference
measurements. The red dotted line shows the performance on GPU without \orcha\ (using OpenACC directives) and the blue dotted line shows the performance on the CPU. Note that both dotted lines represent one measurement only because there is no \datapacket\ generation involved. The solvers operate on one block at a time. For the green line the horizontal axis represents the number of blocks included in one \datapacket.
 
All measurements on GPU used 4 MPI ranks with one GPU per rank.
The GPU performance without using \orcha\ is strikingly poor because a typical block in AMR is 
$16^3$, and when each block ships its data individually 
there is not enough available parallelism to use the GPU
effectively. \orcha\ gets around this problem by combining multiple blocks into one \datapacket\ with 
\milhoja harnessing GPU streams to hide the 
data transfer 
latency.
From the figure one can see that the performance improvement saturates quickly as soon as a \datapacket\ has enough available parallelism.

\cref{tab:sedov} summarizes the measurements of the portions of
the Hydrodynamics solver that can be offloaded to the GPU. The first  mode uses CPU only, and it is the only one that does not use \milhoja. The second mode uses GPU only, while the last three split blocks between the CPU and the GPU in various combinations. The mode names in the table encode how the blocks are split between the two, that is  CPU-$X$ GPU-$Y$ implies
that $X$ blocks are computed on the CPU while $Y$ blocks are shipped
to the GPU. From the table, it is clear 
that for this application on this hardware CPU-only is the best option. 

\begin{table}
    \caption{Time spent in computing the portion of Hydro in the Sedov blast problem in three different modes:
      CPU-only (row 1), GPU-only (row 2), and CPU-GPU-split (rows 3-5).
      In rows 3-5, the split indicates the number of blocks computed simultaneously on CPU and GPU in one cycle.
    All the CPU computations are distributed over 64 cores,
    and the GPU computations are distributed over 4 GPUs in one node of Perlmutter.
    }
    \label{tab:sedov}
    \centering
    \begin{tabular}{cc}
        \toprule
        Mode & Time (seconds) \\
        \midrule
       CPU-only & 7.8  \\
         GPU-only & 13.8 \\
         CPU-10 GPU-20 & 14.3 \\
         CPU-30 GPU-20 & 10.4 \\
         CPU-60 GPU-50 & 9.37 \\
       %\midrule
       \bottomrule
    \end{tabular}
\end{table}

\subsection{Cellular Results}\label{subsec:cullular-results}

\begin{figure}[hb!]
    \centering
    \includegraphics[width=0.9\linewidth]{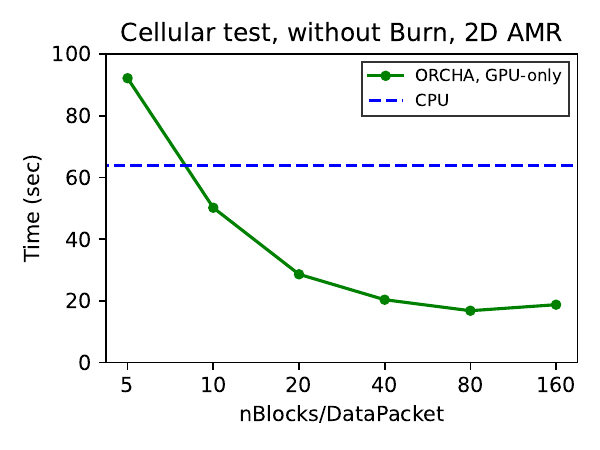}
    \caption{The time spent on Hydrodynamics calculations of
        2D Cellular tests with a varying number of blocks per \datapacket (from 5 to 160 blocks along $x$-axis in log-scale).
        The nuclear burn calculation is turned off to measure
        performance of hydrodynamics on GPU, combined with a complex EOS model
        described as in \cref{sec:cellular}. The higher computational intensity,
        arising from the Helmholtz EOS model, makes the hydrodynamics calculation more
        GPU-favorable, outperforming CPU-only calculation using 64 MPI ranks.
    }
    \label{fig:cellular2d-noburn-perf}
\end{figure}

\cref{fig:cellular2d-noburn-perf} shows the performance of the Cellular application
configured in 2D with Burn turned off to evaluate the
impact of the additional computation introduced by H-EOS. As expected,
the additional computation burden makes GPU more favorable.  \orcha with GPU-only mode
outperforms CPU-only reference. In this figure, we can also more
clearly see the impact of increasing the number of blocks in a
\datapacket. The performance improves until saturation is reached at 80 blocks per \datapacket.

\begin{figure}
    \centering
    \includegraphics[width=0.9\linewidth]{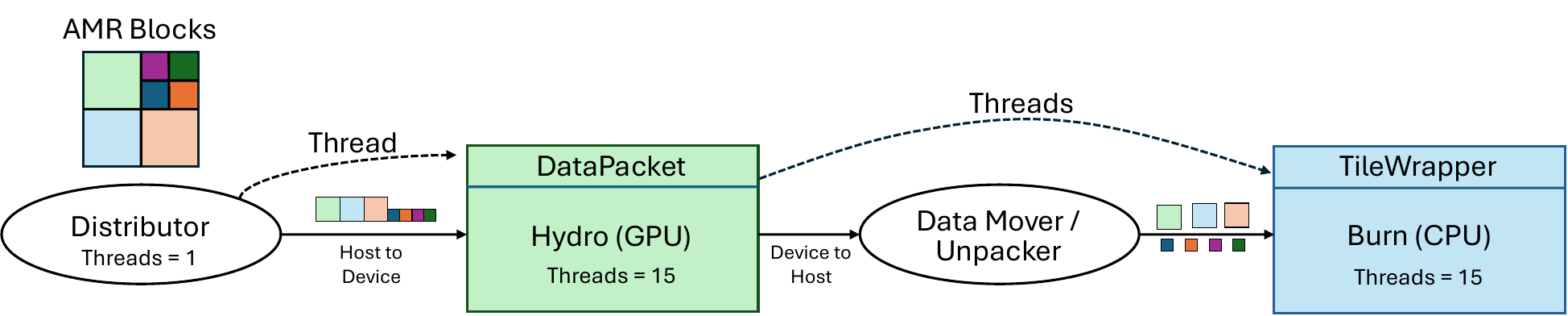}
    \caption{
        \textbf{\extgpu} runtime configuration for Cellular test, calculating Hydro on the GPU,
        followed by running Burn on the CPU.
        The distributor collects several blocks/tiles to consolidate them into a \datapacket,
        then transfers them to the GPU for Hydro calculation. After performing Hydro on the \datapacket,
        the data mover takes the \datapacket, reinterpret the data into \tilewrapper\ form, feed into the
        Burn \taskfunction. The distributor and two \taskfunction{s} can pass their idling threads
        to one another, as described in dotted lines in the figure, to maximally utilize available threads.
        The data dependencies are represented as solid arrows.
    }\label{fig:cellular-extgpu}
\end{figure}

\begin{figure}
    \centering
    \includegraphics[width=0.9\linewidth]{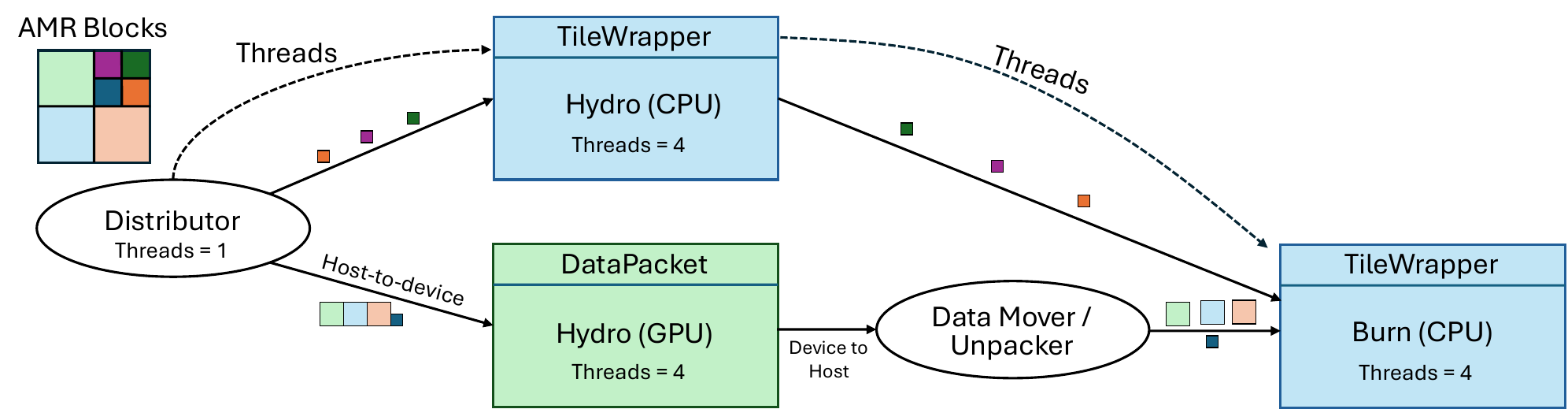}
    \caption{
        \textbf{\extcpugpusplit} runtime configuration for Cellular test,
        calculating Hydro both on the CPU and GPU by splitting AMR blocks, 
        followed by running Burn on the CPU.
        The distributor splits the AMR blocks into a user-defined splitting ratio
        and performs host-to-device transfer. Note that the depicted AMR blocks along the arrows
        are selected for visualizing purposes.
        The solid and dotted lines represent data dependencies and the passing of idling threads, respectively.
    }\label{fig:cellular-extcpugpusplit}
\end{figure}

\begin{figure}
    \centering
    \includegraphics[width=0.9\linewidth]{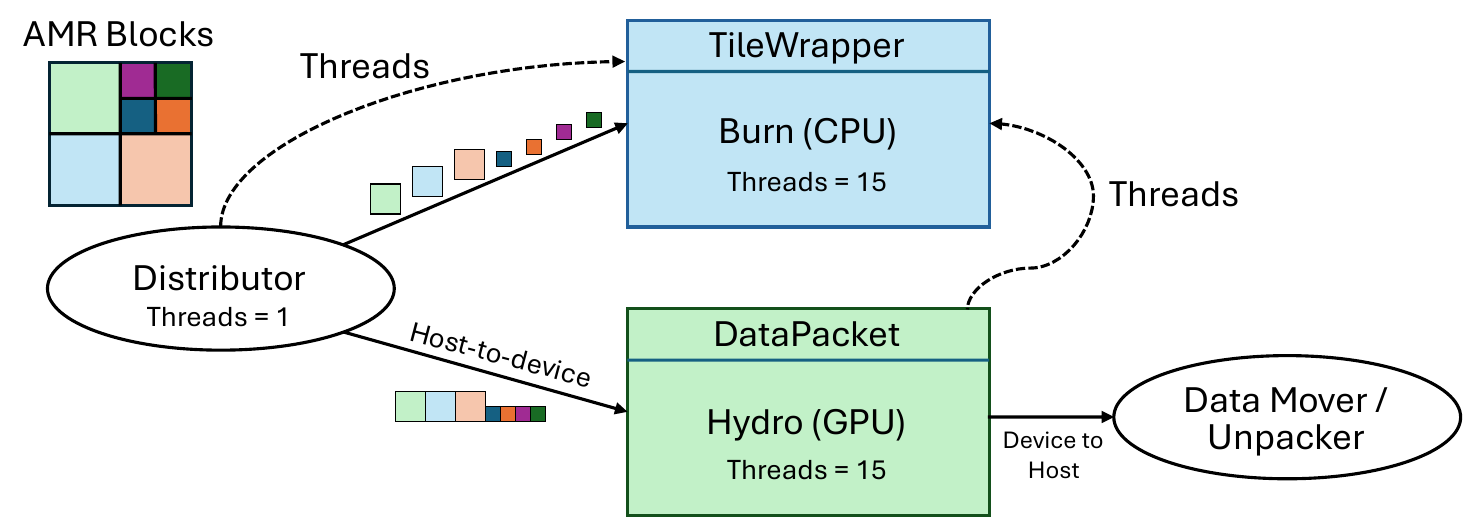}
    \caption{
        \textbf{\cpugpu} runtime configuration for Cellular test,
        calculating Hydro on the GPU while running Burn on the CPU, simultaneously.
        The distributor duplicates AMR blocks and sends them into two \taskfunction{s}
        to exploit the concurrency of the computations.
        Note that even if an AMR block is duplicated for two \taskfunction{s},
        the actual physical data included in the \datapacket\ and \tilewrapper\ may differ.
        For example, pressure information in a block is required for Hydro calculation,
        but it is not needed for Burn, so it will be omitted in its \tilewrapper.
        The solid and dotted lines represent data dependencies and the passing of idling threads, respectively.
    }\label{fig:cellular-cpugpu}
\end{figure}

Other possible \orcha\ enabled   configurations and their corresponding computational pipelines for the Cellular application are enumerated below:
\begin{enumerate}
\item \textbf{\extgpu :} Run Hydro on the GPU, followed by running Burn on the CPU, shown 
  in \cref{fig:cellular-extgpu}.
\item \textbf{\extcpugpusplit :} Run Hydro by splitting AMR blocks between CPU and GPU,
  followed by running Burn on CPU, shown in \cref{fig:cellular-extcpugpusplit}.
\item \textbf{\cpugpu :} Run Hydro on the GPU while running Burn on the CPU simultaneously, shown in \cref{fig:cellular-cpugpu}. This option is less accurate than the first two 
  options because it exercises the cost-benefit trade-off described in \cref{sec:cellular}.
\end{enumerate}

Note that all the above configurations above are obtained by simply changing
the recipe without any modification to the maintained source code in \flashx.
The actual recipes used for configuration \extgpu and \extcpugpusplit are shown in \cref{fig:cellular-recipes}.

\begin{table}[pt!]
    \caption{Time spent in computing Hydro and Burn in the Cellular problem. 
    }\label{tab:cell}
    \centering
    \begin{tabular}{cc}
        \toprule
        mode & Time (seconds) \\
        \midrule
        CPU reference  & 73.779  \\
        \extgpu & 48.874 \\
        \extcpugpusplit & 80.618 \\
        \cpugpu & 45.747 \\
        \bottomrule
    \end{tabular}
\end{table}

 \cref{tab:cell} summarizes the performance of various Cellular  configurations. The first row serves as the reference, listing CPU-only performance without \orcha. Rows two to four list the performances of configurations in \cref{fig:cellular-extgpu,fig:cellular-extcpugpusplit,fig:cellular-cpugpu} respectively. Unlike Sedov, most configurations of Cellular using \orcha\ outperform the reference CPU only mode. In particular, \extgpu\ mode improves performance by approximately 33.8\% over the reference mode.
Conversely, the \extcpugpusplit\ mode shows a performance degradation of
approximately 9.3\% compared to the reference mode.  
While somewhat unexpected,  this outcome also highlights the need to
put a tool in the hands of users that lets them explore different
configurations in planning for production. In modern machines, many
factors come into play during execution, and the obtained performance
does not always match with expectations. In this instance, the
performance degradation is likely to be either because of overheads
from \taskfunction\ synchronization or suboptimal workload
distribution. Knowing about the performance degradation gives users
the choice of either investigating and eliminating the cause or simply
using another configuration that gives a better overall outcome.  

The most notable enhancement was seen in the \cpugpu\ mode,
which concurrently computes Hydro on the GPU while processing Burn on the CPU with a performance gain of about 38.0\%,
demonstrating the effectiveness of parallel execution
and the capability of \orcha\ to leverage heterogeneous hardware resources efficiently.
Overall, the results of our experiments 
demonstrate that \orcha\ can be a very powerful tool to combat the performance portability challenges posed by the dual challenge of simultaneously increasing heterogeneity in hardware and complexity in the scientific applications.  

\begin{figure*}[t]
\centering
% Local override for lstset
{
\lstset{
  basicstyle=\small\ttfamily,
  breaklines=true,
  breakatwhitespace=false,
  columns=fullflexible,
  keepspaces=true,
  breakautoindent=true,
  breakindent=1em,
  postbreak=\mbox{\textcolor{gray}{$\hookrightarrow$}\space},
  numbers=left,
  numbersep=4pt,
  escapeinside={§}{§}
}
\begin{minipage}[t]{0.47\textwidth}
\centering
\caption*{\textbf{\extgpu}}
\vspace{-1em}
\begin{lstlisting}
recipe = flashx.TimeStepRecipe()
begin = recipe.begin_orchestration(
    after=recipe.root,
)

hydro_prepBlock = recipe.add_work(
    "Hydro_prepBlock",
    after=begin,
    map_to=§\codehl{"GPU"}§,
)
hydro_advance = recipe.add_work(
    "Hydro_advance",
    after=hydro_prepBlock,
    map_to=§\codehl{"GPU"}§,
)
burn_burner = recipe.add_work(
    "Burn_burner",
    after=hydro_advance,
    map_to="CPU",
)
burn_update = recipe.add_work(
    "Burn_update",
    after=burn_burner,
    map_to="CPU",
)

end = recipe.end_orchestration(
    begin_node=begin,
    after=burn_update,
)
\end{lstlisting}
\end{minipage}\hfill
\begin{minipage}[t]{0.47\textwidth}
\centering
\caption*{\textbf{\extcpugpusplit}}
\vspace{-1em}
\begin{lstlisting}
recipe = flashx.TimeStepRecipe()
begin = recipe.begin_orchestration(
    after=recipe.root,
)

hydro_prepBlock = recipe.add_work(
    "Hydro_prepBlock",
    after=begin,
    map_to=§\codehl{"CPU,GPU"}§,
)
hydro_advance = recipe.add_work(
    "Hydro_advance",
    after=hydro_prepBlock,
    map_to=§\codehl{"CPU,GPU"}§,
)
burn_burner = recipe.add_work(
    "Burn_burner",
    after=hydro_advance,
    map_to="CPU",
)
burn_update = recipe.add_work(
    "Burn_update",
    after=burn_burner,
    map_to="CPU",
)

end = recipe.end_orchestration(
    begin_node=begin,
    after=burn_update,
)
\end{lstlisting}
\end{minipage}\hfill

} % end local \lstset group
\vspace{-1em}
\caption{
    A side-by-side view of two recipes that were used in testing
    two different hardware configurations for the Cellular detonation test
    in \cref{subsec:cullular-results}. Each physics unit has two API subroutines,
    prefixed as \texttt{Hydro\_} and \texttt{Burn\_}, which are added as nodes in the recipes.
    The corresponding runtime configuration is depicted
    in \cref{fig:cellular-extgpu} and \cref{fig:cellular-extcpugpusplit}, for the recipe on the left and right, respectively.
    The differences between the two recipes are highlighted in yellow.
    Merely minimal recipe changes accomplish two distinct hardware mappings of the exact simulation.
}\label{fig:cellular-recipes}
\end{figure*}

\section{Conclusions}
\label{sec:conclusions}
\begin{figure}
    \centering
    \includegraphics[width=0.9\linewidth]{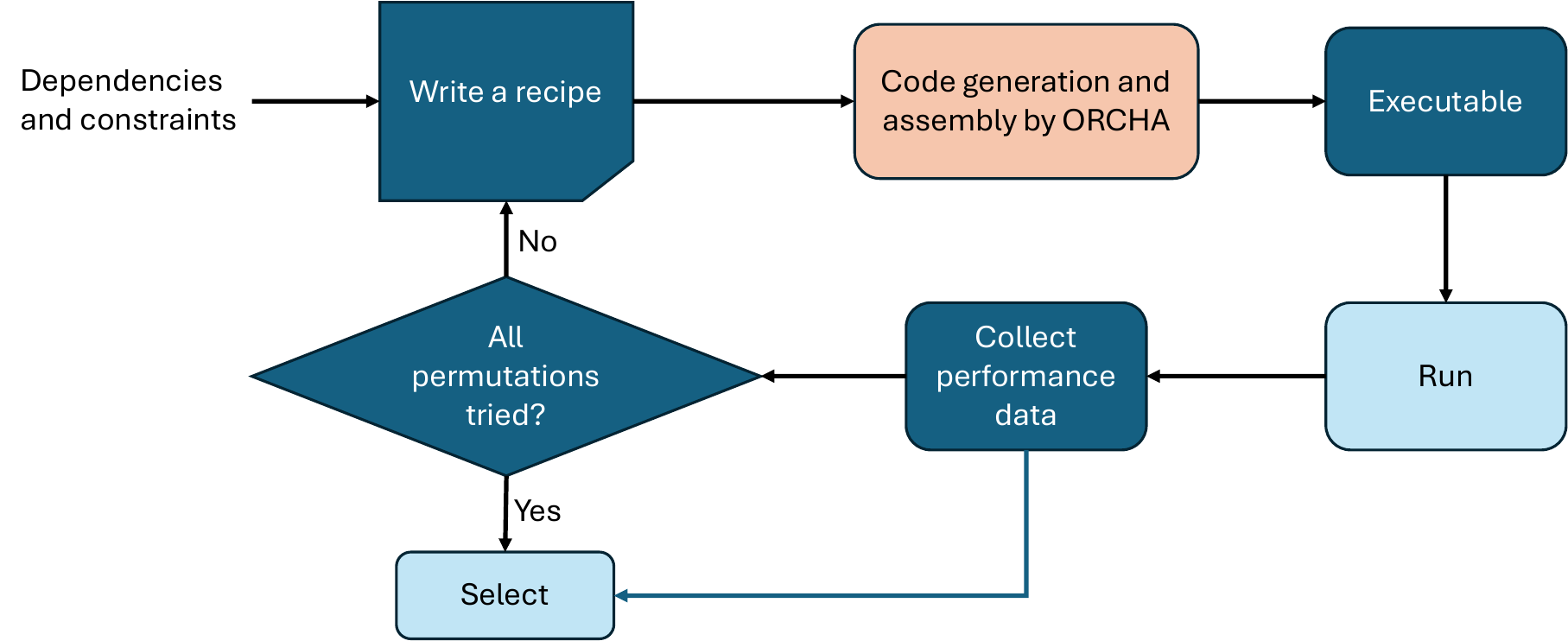}
    \caption{End-to-end workflow for exploring the optimal configuration for a simulation on a specific target.}
    \label{fig:opt-workflow}
\end{figure}

In this paper, we have presented a novel design philosophy for tackling platform heterogeneity that permits a knowledgeable user to explore the best mapping of their computations on heterogeneous platforms.
The key is to hide the nit-picky details through code generation but let the broad strokes of modification and customization remain accessible.
Physics solvers in the code have to follow certain coding standards, such as avoiding dynamic allocations, fully enumerating all the data
in the form of annotations, and ensuring thread safety to be compatible with \orcha. In return, they get the option of exploring the best possible configuration for their science runs by
simply rewriting the recipes without the necessity of modifying the maintained source code as shown in \cref{fig:opt-workflow}. The degree of customizability is controlled by the authors of the physics solver through the interfaces they expose to \orcha. 

Individual tools that are part of \orcha are not necessarily optimal. For example, there is ample scope for kernel-level refinements
or advanced runtime tuning. However, the design of \orcha ensures that these advances can be made in the tools individually and gradually, or even swap out or add more tools in the end-to-end pipeline.
Our initial tests already yield modest, repeatable speedups when GPU configurations are generated from the recipe and compared against a CPU-only baseline. They also yield performance gain on the GPU by consolidating several AMR blocks into one \datapacket over using OpenMP/OpenACC offloading on a block-by-block basis. 
\orcha is designed to grow more sophisticated by adding different tools with different capabilities by letting the tools remain composable. Therefore, our future work will focus on integrating tools that can do auto-tuned kernel generation, enhance memory layout, and adaptive scheduling policies into the runtime. 

\section*{Acknowledgment}
The submitted manuscript has been created by UChicago Argonne, LLC,
Operator of Argonne National Laboratory (“Argonne”). Argonne, a
U.S. Department of Energy Office of Science laboratory, is operated
under Contract No. DE-AC02-06CH11357. The U.S. Government retains for
itself, and others acting on its behalf, a paid-up nonexclusive,
irrevocable worldwide license in said article to reproduce, prepare
derivative works, distribute copies to the public, and perform
publicly and display publicly, by or on behalf of the Government.  The
Department of Energy will provide public access to these results of
federally sponsored research in accordance with the DOE Public Access
Plan. http://energy.gov/downloads/doe-public-access-plan. 
This work was supported by the Exascale Computing Project
(17-SC-20-SC), a collaborative effort of the U.S.  Department of
Energy Office of Science and the National Nuclear Security
Administration, and by the Scientific Discovery through Advanced
Computing (SciDAC) program via the  Office of Nuclear Physics and
Office of Advanced Scientific Computing Research in the Office of
Science at the U.S.\ Department of Energy. Also, this research used
resources of the National Energy Research Scientific Computing Center,
which is supported by the Office of Science of the U.S. Department of
Energy under Contract No. DE-AC02-05CH11231.

\noindent{\bf Funding} The author(s) disclosed receipt of the
following financial support for the research, authorship, and/or
publication of this article: U.S. Department of Energy, Office of
Science, under contract number DE-AC02-06CH11357

\bibliographystyle{SageH}
% \bibliography{references.bib}

% NOTE: ScholarOne system doesn't work well with the .bib file,
%       although their LaTeX template file uses BibTeX.
%       I copied the contents of the generated .bbl file here
%       instead of using the BibTeX.

\clearpage

\appendix

\section{Appendix}

\subsection{Example Macro Usage in \orcha}\label{appendix:macro}
This appendix is dedicated to describing the macro constructs used in \orcha,
with a focus on the syntax introduced in \cref{lst:contrived_defs}.
The examples provided are intentionally simplified and contrived,
and they do not represent the actual application code.
Instead, they are designed to illustrate the use and expansion of macros
that can be used in an application with ORCHA.

Before we start the discussion, there are two points to clarify. First, all macro names and their definitions shown below are chosen by the application developer and are not built into \macro. Other applications with a different data structure may have distinct sets of macros with different definitions. Second, \macro only performs textual substitution of macro identifiers and handles arbitrations among alternative definitions; it has no knowledge of the target language or code logic. The example shown below was achieved entirely by providing two alternative macro definitions for CPU and GPU, not by any built-in hardware awareness of \macro.

In this example, we consider a subroutine \texttt{P1}, illustrated in \cref{lst:contrived_defs}:
\begin{lstlisting}[
    language=Fortran,
    caption={Example source code defined with macros},
    label={lst:appendix_p1_defs}
]
subroutine P1(Ustate, Udesc, Scr1, Scr2, S1, S2, S3)
    @M Udesc_declare
    real, dimension(@M Udesc_view), intent(inout) :: Ustate, Scr1, Scr2
    real, intent(in) :: S1, S2, S3
    integer :: i, j, k
    @M loop_ijk_begin
       @M Comp1(Ustate, Scr1, Scr2, S1)
    @M loop_ijk_spl_end
    @M loop_ijk_spl_begin
       @M Comp2(Ustate, Scr1, S2, S3)
    @M loop_ijk_end
end subroutine P1
\end{lstlisting}
The macros are encoded with \texttt{@M} keyword, followed by their name.
The \texttt{@M} macro keyword is a design choice in this example, and can be readily modified within the \macro if desired.
Since \orcha can launch this API subroutine on either CPU or GPU, determined by a user-provided recipe,
we are interested in generating both CPU and GPU versions of the subroutine,
which potentially have different code structures for the same computations.

Macro definitions are specified in separate files and are passed to \macro with the selected macro-annotated source files to expand.
Macro definitions are specified using a structured form as in \cref{lst:appendix_macro_defs_format}, where it defines a macro by its identifier, optional arguments, and associated expansion statements.
\begin{lstlisting}[
    language={},
    caption={Macro definition format},
    label={lst:appendix_macro_defs_format}
]
[<macro name>]
(optional) args = <arg1>, <arg2>, ...
definitions =
  <expansion statements>
\end{lstlisting}

We first define the macros that are used both for CPU and GPU, with the same definitions:
\begin{lstlisting}[
    language={},
    caption={Common macros for \texttt{P1}.},
    label={lst:appendix_common_macros}
]
[Udesc_declare]
definition =
  type(Udesc_type), intent(in) :: Udesc

[Udesc_view]
definition =
  Udesc%lo(1):Udesc%hi(1), Udesc%lo(2):Udesc%hi(2), Udesc%lo(3):Udesc%hi(3)
\end{lstlisting}
The two macros above show an example of macro usage to minimize error-prone repetition of common code patterns. 

We now have two distinct sets of macro definitions that alter \texttt{P1} for the CPU and GPU variants. For CPU, for instance, we have:
\begin{lstlisting}[
    language={},
    caption={Macro definitions for CPU version of \texttt{P1}.},
    label={lst:appendix_cpu_macros}
]
[loop_ijk_begin]
definition =
  do k = Udesc%lo(3), Udesc%hi(3)
     do j = Udesc%lo(2), Udesc%hi(2)
        do i = Udesc%lo(1), Udesc%hi(1)

[loop_ijk_spl_end]
definition =

[loop_ijk_spl_begin]
definition =

[loop_ijk_end]
definition =
        end do
     end do
  end do

[Comp1]
args = U, scr1, scr2, s
definition =
  call Comp1_CPU(U, scr1, scr2, s, i, j, k)

[Comp2]
args = U, scr, s1, s2
definition =
  call Comp2_CPU(U, scr, s1, s2, i, j, k)
\end{lstlisting}
Note that the macros \texttt{loop\_ijk\_spl\_end} and \texttt{loop\_ijk\_spl\_begin} have null definitions.
This indicates that \macro will not expand those two macros, ignoring them in expansion \texttt{P1} for the CPU version.
The macros \texttt{Comp1} and \texttt{Comp2} accept arguments, which are substituted with the provided actual arguments during the expansions.

On the contrary, we have a different set of macro definitions for the GPU version of \texttt{P1}.
\begin{lstlisting}[
    language={},
    caption={Macro definitions for GPU version of \texttt{P1}.},
    label={lst:appendix_gpu_macros}
]
[loop_ijk_begin]
definition =
  !$acc parallel loop collapse(3) present(Ustate, Scr1, Scr2)
  do k = Udesc%lo(3), Udesc%hi(3)
     do j = Udesc%lo(2), Udesc%hi(2)
        do i = Udesc%lo(1), Udesc%hi(1)

[loop_ijk_spl_end]
definition =
  @M loop_ijk_end

[loop_ijk_spl_begin]
definition =
  !$acc parallel loop collapse(3) present(Ustate, Scr1)
  do k = Udesc%lo(3), Udesc%hi(3)
     do j = Udesc%lo(2), Udesc%hi(2)
        do i = Udesc%lo(1), Udesc%hi(1)

[loop_ijk_end]
definition =
        end do
     end do
  end do

[Comp1]
args = U, scr1, scr2, s
definition =
  call Comp1_GPU(U, scr1, scr2, s, i, j, k)

[Comp2]
args = U, scr, s1, s2
definition =
  call Comp2_GPU(U, scr, s1, s2, i, j, k)
\end{lstlisting}
The notable difference in the GPU macro definitions is that they provide explicit definitions
for \texttt{loop\_ijk\_spl\_end} and \texttt{loop\_ijk\_spl\_begin}, to split the loop into two,
so that two distinct computations happen in two different kernels. Note that the macro \texttt{loop\_ijk\_spl\_end}
is defined in terms of another macro, which is recursively expanded during macro expansion by \macro;
this is allowed as long as no cyclic definitions are introduced.

Finally, two variants of \texttt{P1}, expanded from \cref{lst:appendix_p1_defs}
with different macro definitions described above, are presented in \cref{fig:appendix_side_by_side}.

\begin{figure*}[t]
\centering
% Local override for lstset
{
\lstset{
  basicstyle=\small\ttfamily,
  breaklines=true,
  breakatwhitespace=false,
  columns=fullflexible,
  keepspaces=true,
  breakautoindent=true,
  breakindent=1em,
  postbreak=\mbox{\textcolor{gray}{$\hookrightarrow$}\space},
  numbers=left,
  numbersep=4pt,
  escapeinside={§}{§}
}
\begin{minipage}[t]{0.47\textwidth}
\centering
\caption*{CPU variant of \texttt{P1}.}
\vspace{-1em}
\begin{lstlisting}[
    language=Fortran
]
subroutine P1_CPU(Ustate, Udesc, Scr1, Scr2, S1, S2, S3)
   type(Udesc_type), intent(in) :: Udesc
   real, dimension(Udesc%lo(1):Udesc%hi(1), Udesc%lo(2):Udesc%hi(2), Udesc%lo(3):Udesc%hi(3)), intent(inout) :: Ustate, Scr1, Scr2
   real, intent(in) :: S1, S2, S3
   integer :: i, j, k

   do k = Udesc%lo(3), Udesc%hi(3)
      do j = Udesc%lo(2), Udesc%hi(2)
         do i = Udesc%lo(1), Udesc%hi(1)
            call Comp1_CPU(Ustate, Scr1, Scr2, s1, i, j, k)
            call Comp2_CPU(Ustate, Scr1, s2, s3, i, j, k)
         end do
      end do
   end do
end subroutine P1
\end{lstlisting}
\end{minipage}\hfill
\begin{minipage}[t]{0.47\textwidth}
\centering
\caption*{GPU variant of \texttt{P1}.}
\vspace{-1em}
\begin{lstlisting}[
    language=Fortran
]
subroutine P1_GPU(Ustate, Udesc, Scr1, Scr2, S1, S2, S3)
   type(Udesc_type), intent(in) :: Udesc
   real, dimension(Udesc%lo(1):Udesc%hi(1), Udesc%lo(2):Udesc%hi(2), Udesc%lo(3):Udesc%hi(3)), intent(inout) :: Ustate, Scr1, Scr2
   real, intent(in) :: S1, S2, S3
   integer :: i, j, k

   !$acc parallel loop collapse(3) present(Ustate, Scr1, Scr2)
   do k = Udesc%lo(3), Udesc%hi(3)
      do j = Udesc%lo(2), Udesc%hi(2)
         do i = Udesc%lo(1), Udesc%hi(1)
            call Comp1_GPU(Ustate, Scr1, Scr2, s1, i, j, k)
         end do
      end do
   end do

   !$acc parallel loop collapse(3) present(Ustate, Scr1)
   do k = Udesc%lo(3), Udesc%hi(3)
      do j = Udesc%lo(2), Udesc%hi(2)
         do i = Udesc%lo(1), Udesc%hi(1)
            call Comp2_CPU(Ustate, Scr1, s2, s3, i, j, k)
         end do
      end do
   end do
end subroutine P1
\end{lstlisting}
\end{minipage}\hfill

} % end local \lstset group
\vspace{-1em}
\caption{
    Two generated variants of \texttt{P1} expanded by \macro with different macro definitions.
}\label{fig:appendix_side_by_side}
\end{figure*}

\end{document}

%% file: preamble.tex
\usepackage{moreverb}

\usepackage{xr}

\usepackage{xurl}
\usepackage[colorlinks,bookmarksopen,bookmarksnumbered,citecolor=red,urlcolor=red,breaklinks=true]{hyperref}

\usepackage{doi}
\usepackage[hyphenbreaks]{breakurl}

\usepackage{graphicx}
\usepackage{color} % To define colors
\usepackage{setspace}
\usepackage{natbib}
\usepackage{subcaption}
\usepackage{xspace}
\usepackage{tabularx}
\usepackage{cite}
\usepackage{amsmath,amssymb,amsfonts}
\usepackage{algorithmic}
\usepackage{graphicx}
\usepackage{textcomp}
\usepackage[table]{xcolor}
\usepackage{soul}
\usepackage{listings}  % this must be included BEFORE cleveref
\usepackage[capitalise]{cleveref}
\usepackage{booktabs}

\def\volumeyear{2025}

\newcommand{\spark}{Spark\xspace}
\newcommand{\flashx}{\mbox{Flash-X}\xspace}
\newcommand{\orcha}{ORCHA\xspace}

\newcommand{\cgkit}{\mbox{CG-Kit}\xspace}
\newcommand{\milhoja}{Milhoja\xspace}
\newcommand{\macro}{Macroprocessor\xspace}
\newcommand{\recipetools}{FlashX-RecipeTools\xspace}

\newcommand{\dataitem}{DataItem\xspace}
\newcommand{\datapacket}{DataPacket\xspace}
\newcommand{\taskfunction}{TaskFunction\xspace}
\newcommand{\tilewrapper}{TileWrapper\xspace}
\newcommand{\extgpu}{$\text{(GPU)}_\text{Hydro} \to \text{(CPU)}_\text{Burn}$}
\newcommand{\cpugpu}{$\text{(GPU)}_\text{Hydro} \parallel \text{(CPU)}_\text{Burn}$}
\newcommand{\extcpugpusplit}{$\text{(GPU+CPU)}_\text{Hydro} \to \text{(CPU)}_\text{Burn}$}

\newcommand{\codehl}[1]{\colorbox{yellow!55}{#1}}